\documentstyle[12pt, verbatim]{article}
\font\teneuf=eufm10 
\font\fiveeuf=eufm5
\font\seveneuf=eufm7
\newfam\euffam
\def\euf{\fam\euffam\teneuf}
\textfont\euffam=\teneuf
\scriptfont\euffam=\seveneuf
\scriptscriptfont\euffam=\fiveeuf 
\renewcommand{\baselinestretch}{1.3}

.tfm scaled \magstep1
.tfm 

\catcode`@=11
\newskip\@centering \@centering=0pt plus 1000pt minus 1000pt
\def\openup{\afterassignment\@penup\dimen@=}
\def\@penup{\advance\lineskip\dimen@
  \advance\baselineskip\dimen@
  \advance\lineskiplimit\dimen@}
\def\eqalign#1{\null\,\vcenter{\openup\jot\m@th
  \ialign{\strut\hfil$\displaystyle{##}$&$\displaystyle{{}##}$\hfil
      \crcr#1\crcr}}\,}
\newif\ifdt@p
\def\displ@y{\global\dt@ptrue\openup\jot\m@th
  \everycr{\noalign{\ifdt@p \global\dt@pfalse
      \vskip-\lineskiplimit \vskip\normallineskiplimit
      \else \penalty\interdisplaylinepenalty \fi}}}
\def\@lign{\tabskip\z@skip\everycr{}} 
\def\displaylines#1{\displ@y
  \halign{\hbox to\displaywidth{$\@lign\hfil\displaystyle##\hfil$}\crcr
    #1\crcr}}
\def\eqalignno#1{\displ@y \tabskip\@centering
  \halign to\displaywidth{\hfil$\@lign\displaystyle{##}$\tabskip\z@skip
    &$\@lign\displaystyle{{}##}$\hfil\tabskip\@centering
    &\llap{$\@lign##$}\tabskip\z@skip\crcr
    #1\crcr}}
\def\leqalignno#1{\displ@y \tabskip\@centering
  \halign to\displaywidth{\hfil$\@lign\displaystyle{##}$\tabskip\z@skip
    &$\@lign\displaystyle{{}##}$\hfil\tabskip\@centering
    &\kern-\displaywidth\rlap{$\@lign##$}\tabskip\displaywidth\crcr
   #1\crcr}}
\catcode`@=12

\def \mod{{\rm mod}\,}

\def\e{\varepsilon}

\begin{document}

\title{\large\bf THE DIMENSION OF THE SPACE OF CUSP FORMS \\
                 OF WEIGHT ONE}

\author{\large \bf W. Duke%
\thanks{Research supported in part by NSF Grants
          DMS--9202022 and DMS--9022140}
}
\maketitle

\subsection*{1. Introduction} 

It is a basic problem to determine the dimension of the space of cusp
forms of a given type. For classical holomorphic forms of integral
weight larger than one, the dimension is well understood by means of
either the Riemann--Roch theorem or the Selberg trace formula.  The
case of weight one, however, remains mysterious. From the point of
view of spectral theory this is because these forms belong to an
eigenvalue of the Laplacian which is not isolated; the difficulty of
estimating nontrivially the multiplicity of such an eigenvalue is well
known.

\medskip
Suppose, for example, that $q$ is a prime and that $S_1(q)$ is the
space of holomorphic cusp forms for $\Gamma_0(q)$ of weight 1 with
character $\left({\cdot\over q}\right )$, the Legendre symbol. No
nonzero cusp forms may exist unless $\left({-1\over q}\right ) = -1$,
so assume $q\equiv 3\; (\mod \, 4)$. Hecke discovered that the
existence of such cusp forms is tied up with the class number $h$ of
${\bf Q}(\sqrt{-q})$; if $\chi$ is any nontrivial (hence non--real)
class character then
$$
\sum_{{\euf a}} \chi({\euf a}) e(N({\euf a})\, z) \in S_1(q)\,,
\leqno(1)
$$
where ${\euf a}$ runs over all nonzero integral ideals of 
${\bf Q}(\sqrt{-q})$. 
There are $(h-1)/2$ independent forms of this type so by Siegel's
theorem we have the (ineffective) lower bound
$$
{\rm dim}\, S_1(q) \mathop{\gg}_{\varepsilon} q^{1/2-\e}
$$
for all $\e > 0$.

In general, $S_1(q)$ is not spanned by forms of Hecke's type. The
construction of specific examples demonstrating this is an active area
of research (see \cite{F} and its references). Such exotic forms seem
quite rare, however, and it appears reasonable to expect that in fact
$$
{\rm dim}\, S_1(q) =\frac{1}{2}\, (h-1)+O(q^\e)\,.
$$
In particular, this would imply that ${\rm dim}\, S_1(q) \ll q^{1/2}
\log q$.  Rather less, however, is actually known. Serre has shown
\cite{S} that if $q = 24m-1$ or $24m + 7$ then ${\rm dim}\, S_1(q) \le
m - (h-1)/2$, while if $q = 24m + 11 $ or $24m + 19$ then ${\rm dim}\,
S_1(q) \le m - h+1$.  Sarnak has informed me that the
Selberg trace formula for weight one \cite[Chapter 9]{Hej} with a
suitably chosen test function yields the bound
$$
{\rm dim}\, S_1(q)\ll \frac{q}{\log q}
$$
which, for large values of $q$, is currently the best
known.\footnote{In an unpublished note J.-M. Deshouillers and
H. Iwaniec obtained in this way the bound $O(q \log^{-3} q)$ for the
multiplicity of the eigenvalue $1/4$ in the case of weight zero Maass
cusp forms for $\Gamma_0(q)$ with trivial character.}  The main object
here is to improve this estimate.

\medskip\noindent
{\bf Theorem 1.} {\sl For $q$ prime
$$
{\rm dim}\, S_1(q)\ll q^{11/12}\, \log^4 q\,,
$$
with an absolute implied constant.
}

\bigskip
Roughly speaking, the idea of the proof is to exploit two conflicting
properties of the Fourier coefficients of newforms in $S_1(q)$ not of
Hecke's type (1): their approximate orthogonality and the finiteness
of the number of their possible values at primes. The first property
is a consequence of their belonging to automorphic forms while the
second is a consequence of Deligne--Serre theorem. Taken together,
these properties limit the number of possible newforms which may
exist.

Following Serre \cite{S}, Theorem~1 has an application to the quotient
$X_0^*(q)$
of the modular curve
$X_0(q)$
by the Fricke involution $z\mapsto -1/qz$ when 
$q = 24 m - 1$ is prime. In this case the genus of 
$X_0^*(q)$
is $m-(h-1)/2$.

\bigskip\noindent
{\bf Corollary.} {\sl For $q = 24m - 1$ prime, the space of differential
forms of the first kind on $X_0^*(q)$ with a zero of order
at least $m$ at the cusp has dimension which is 
$O(q^{11/12} \log^4 q)$.
}

\bigskip
A different kind of application of these ideas is to bound the number
$m _4(-q)$ of quartic number fields of discriminant $-q$.

\bigskip\noindent
{\bf Theorem 2.} {\sl For $q$ prime
$$
m _4(-q) \ll q^{7/8} \log ^4 q\,,
$$
with an absolute implied constant.
}

\bigskip\noindent
To put this result in context it may be of interest to show what follows
from algebraic number theory and trivial bounds for class numbers.
By means of class field theory Heilbronn \cite{Hei} showed that 
$$
m_4(-q) = \frac{4}{3}\, \sum_k h_2(k)\,,
$$
where $k$ runs over all cubic number fields of discriminant $-q$.
Here, for any $\ell$ and any number field $k$, $h_\ell(k)$ denotes the
number of ideal classes of $k$ of (exact) order $\ell$. Furthermore, by
\cite{Ha}, the number of cubic fields in the sum is $\frac{3}{2}\,
h_3({\bf Q}(\sqrt{-q}))$ . For the class number $h(k)$ of any number
field $k$ of degree $n>1$ and discriminant $D$ we have the bound
$$
h(k) \ll |D|^{1/2}\log^{n-1}|D|
$$
where the implied constant depends only on $n$ (see
\cite[p.153.]{N}).  Since $h_\ell(k) \leq h(k)$ we deduce that
$$
m _4(-q) \ll q\, \log^3 q 
$$
with an absolute implied constant.  The improvement of this given in
Theorem~2 requires both the classification of quartic fields of
discriminant $-q$ by odd octahedral Galois representations of
conductor $q$ given in \cite{S} and the proof in this case of the
Artin conjecture given in \cite{T}.  If we assume the Artin conjecture
for icosahedral representations then similarly we can deduce that the
number of non-real quintic fields of discriminant $q^2$ whose normal
closure has Galois group $A_5$ is $O(q^{11/12}\log^4 q)$.

\subsection*{2. Approximate orthogonality of Fourier Coefficients}

For $q\in {\bf Z}^+$ and $\e$ an odd Dirichlet character $\mod q$,
let $S_1(q)=S_1(q,\e)$ be the set of all holomorphic cusp forms for
$\Gamma_0(q)$ of weight 1 with character $\e$. Thus $f\in S_1(q)$
satisfies for $\gamma\in {a\; b\choose c\; d}\in \Gamma_0(q)$
$$
f(\gamma z) = \e(d) \, (cz+d)\, f(z)
$$
and $({\rm Im} \, z)^{1/2} \, |f(z)|$ is uniformly bounded on the 
upper half--plane ${\bf H}$. The vector space $S_1(q)$ is finite 
dimensional and has an inner product
$$
\langle f, g\rangle 
=
\int\limits_{\Gamma_0(q)\backslash {\bf H}}
f(z) \bar g(z)\, y^{-1}dxdy\,.
$$
Each $f\in S_1(q)$ has the Fourier expansion at $\infty$
$$
f(z) =\sum_{n\ge 1} a_f(n)\, e(nz)\,.
$$
The object of this section is to establish the following mean value
result which expresses the approximate orthogonality of the $a_f(n)$
over any fixed orthonormal basis ${\cal B}$ for $S_1(q)$. 

\medskip\noindent
{\bf Proposition 1.} {\sl For arbitrary $c_n\in {\bf C}$ with 
$1\le n\le N$ we have
$$
\sum_{f\in {\cal B}} \, \left|\sum_{n\le N} c_n\, a_f(n)\right|^2
\ll
\left(1+\frac{N}{q}\right)\, \sum_{n\le N}\, |c_n|^2
$$
with an absolute implied constant.
}

\medskip\noindent
The proof we give of this uses the following duality principle.

\medskip\noindent
{\bf Lemma 1.} {\sl Suppose that $V$ is a finite dimensional inner 
product space over ${\bf C}$ with an orthonormal basis ${\cal B}$.
Let $\{v_1,\ldots,v_N\}\subset V$ be a fixed set of vectors and $\Delta$
be a positive number.  Then the inequality
$$
\sum_{n\le N}\, |\langle u,v_n \rangle |^2 \le \Delta\, \langle u,u\rangle
\leqno(2)
$$
holds for all $u\in V$ if and only if 
$$
\sum_{f\in {\cal B}}\, \left| \sum_{n\le N} c_n\, \langle f, v_n\rangle
 \right|^2\le \Delta\, \sum_{n\le N}\, |c_n|^2
\leqno(3)
$$
holds for all $c_n\in {\bf C}$.
}

\medskip\noindent 
{\bf Proof:} By Parseval's equality
$$
\eqalign{
\sum_{f\in {\cal B}} \left|\sum_{n \le N} c_n \, \langle f, v_n\rangle
 \right|^2
&=
\mathop{\sum\sum}_{n,m\le N}\, c_n\bar c_m
\sum_{f\in {\cal B}} \langle f, v_n\rangle\, 
\overline{\langle f, v_m\rangle}
\cr\noalign{\vskip 5pt}
&=
\mathop{\sum\sum}_{n,m\le N} c_n \bar c_m \langle v_n,v_m\rangle\,.
\cr
}
\leqno(4)
$$
Now the Lemma follows from Boas' generalization of Bessel's inequality
as given in \cite[p.~151]{D}.  For the sake of completeness we will
include the rest of the proof that $(2)\Rightarrow (3)$ here since
this is the part which we need and refer to \cite{D} for the converse.

Take $u=\sum_{m\le N} c_m v_m$ in (2); then
$$
\eqalign{
\langle u, u\rangle &=
\mathop{\sum\sum}_{n,m\le N} c_n \bar c_m\langle v_n, v_m\rangle
=\sum_{n\le N} c_n \langle v_n, u\rangle 
\cr\noalign{\vskip 5pt}
&\le
\left(\sum_{n\le N} |\langle u, v_n\rangle|^2\right)^{1/2}\, 
\left(\sum_{n\le N}|c_n|^2\right)^{1/2}
\cr\noalign{\vskip 5pt}
&\le
\Delta^{1/2} \langle u, u\rangle ^{1/2}\, \left(\sum |c_n|^2\right)^{1/2}
\cr
}
\leqno(5)
$$
by (2). Hence
$$
\langle u, u\rangle \le \Delta\, \sum_{n\le N} |c_n |^2
$$
and (3) follows by (4) and (5). $\Box$

\medskip
Taking $v_n = \sum_{f\in {\cal B}} \bar a_f(n)\, f$, we see that
Proposition~1 is reduced to the following Lemma. 

\medskip\noindent
{\samepage {\bf Lemma 2.} {\sl For any $f\in S_1(q)$
$$
\sum_{n \le N} \, |a_f(n)|^2
\ll 
\left(1 +\frac{N}{q}\right)\, \langle f, f\rangle
$$
with an absolute implied constant.
}}

\medskip\noindent
{\bf Proof:} We employ a technique of Iwaniec (unpublished)
to bound such sums uniformly. For any $y> 0$ we have
$$
\sum_{n\ge 1} e^{-4\pi n y} \, |a_f(n)|^2
=\int_0^1 |f(x + iy)|^2 \, dx.
$$
Thus, for any $Y>0$,
$$
\sum_{n\le N} F_Y(n)\, |a_f(n)|^2 \leq \int _Y^\infty \int _0^1
|f(z)|^2 y^{-1} \, dxdy\,,
\leqno(6)
$$
where
$$
F_Y(n)=\int_1^\infty e^{-4\pi n Y y}y^{-1} \, dy\,
\ge \int_1^\infty e^{-4\pi N Yy}\, y^{-1}\, dy
\leqno(7)
$$
for $n\le N$. Setting 
$$P(Y) = \{z\in {\bf H}\, ;\, 0 < {\rm Re}\,
z\le 1,\; {\rm Im}\, z> Y\}
$$ 
we have
$$
\int_Y^\infty \int_0^1 |f(z)|^2\, y^{-1}\, dxdy
\le \mathop{\max}_{z\in P(Y)\,}\# \{\gamma \in \Gamma_0(q)/\{\pm 1\}\, ;\, \gamma z\in P(Y)\}
\, \langle f, f\rangle \,.
\leqno(8)
$$
Now for fixed $z=x+iy \in {\bf H}$ the condition that $\gamma z = (az+b)/(cqz+d)
\in P(Y)$ is imposed by requiring that
$$
{\rm Im}\,((az+b)/(cqz+d)) = y |cqz+d|^{-2}>Y
\leqno(9)
$$
and that 
$$
0<{\rm Re}\,((az+b)/(cqz+d))\leq 1.
\leqno(10)
$$
For a fixed $c$ and $d$
satisfying (9) $a$ and $b$ are determined by (10).
It follows that for $z \in {\bf H}$
$$
\eqalign{
\#\{\gamma \in \Gamma_0(q)/\{\pm 1\}\, ; \gamma z\in P(Y)\} &\le 1 + \# \{c > 0, d\,;\, |cqz+d|^2 < yY^{-1}\}
\cr\noalign{\vskip 5pt}
&\le 1 +\# \{c, d\,;\, 0<c< q^{-1}(yY)^{-1/2}\; {\rm and}\; |cqx+d|^2
    < yY^{-1}\} }
$$
by taking real and imaginary parts. Hence for $z \in P(Y)$
$$
\eqalign{
\#\{\gamma \in \Gamma_0(q)/\{\pm 1\}\}\, ; \gamma z \in P(Y)\} &\le 1 + q^{-1}(yY)^{-1/2}\,
    (1 + 2y^{1/2} Y^{-1/2}) \cr\noalign{\vskip 5pt} &= 1 +
    q^{-1}(yY)^{-1/2} + 2 q^{-1} Y^{-1} \cr\noalign{\vskip 5pt} &\le 1
    + 3 q^{-1} Y^{-1}.}
$$
Thus from (6) and (8) we get
$$
\sum_{n\le N} F_Y(n)\, |a_f(n)|^2
\le 
(1+3q^{-1} Y^{-1})\, \langle f,f\rangle\,.
$$ 
Choosing $Y=3N^{-1}$ and using (7) we get Lemma~2 (with $10^{18}$ for
the absolute constant). $\Box$

\subsection*{3. Consequences of the Deligne--Serre theorem}

Let $N_1(q, \e)\subset S_1(q,\e)$ be the set of normalized newforms.
For $f\in N_1(q, \e)$ the associated $L$--function is an Euler product
$$
L_f(s) 
=
\sum_{n\ge 1} a_f(n)\, n^{-s} 
=
\prod_p \, (1-a_f(p)\, p^{-s}+\e(p)\, p^{-2s})^{-1}\,.
\leqno(11)
$$
The Deligne--Serre theorem \cite{DS} states that $L_f(s)$ is the 
Artin $L$--function of an irreducible two--dimensional Galois 
representation $\rho$ of conductor $q$ with $\det \rho=\e$ (via 
the Artin map). We shall use two consequences of this result.
The first is that $a_f(n)$ satisfies the Ramanujan bound
$$
|a_f(n)|\le d(n),
\leqno(12)
$$
where $d(n)$ is the divisor function.
The second is that $f$ may be classified as being of dihedral,
tetrahedral, octahedral, or icosahedral type according to whether the
image of $\rho$ in $PGL(2, {\bf C})$ is $D_h$, $A_4$, $S_4$ or $A_5$.

We now restrict attention to the case that $q$ is prime and
$\e(\cdot)=\left(\frac{\cdot}{q}\right)$ and write $N_1(q)$ for
$N_1(q, \e)$. Also, let $N_{{\rm dih}}, N_{{\rm oct}}, N_{{\rm ico}}$
be the forms in $N_1(q)$ of dihedral, octahedral and icosahedral type,
respectively. It is shown in \cite[p.~343 of Collected Papers]{S} that
$f$ is of dihedral type exactly when it is of Hecke's type (1).

\medskip\noindent
{\bf Proposition 2.} {\sl Suppose that $f\in N_1(q)$ for $q$
prime. Then

\smallskip\noindent
(a)
$\langle f, f\rangle \ll q\log^3 q$ with an absolute implied constant.

\smallskip\noindent
(b) If $f$ is not of dihedral type then it is either of octahedral or
icosahedral type.  If $f\in N_{{\rm oct}}$ then
$$
a_f(p^8)-a_f(p^4)-\left(\frac{p}{q}\right)\, a_f(p^2)=1\,,
\leqno(13)
$$
while if $f\in N_{{\rm ico}}$ then
$$
a(p^{12})-a_f(p^8)-\left(\frac{p}{q}\right)\, a_f(p^2)=1
\leqno(14)
$$
for all primes $p\ne q$.
}

\medskip\noindent
{\bf Proof:} (a) The Rankin--Selberg convolution
$$
\varphi(s)=\sum_{n\ge 1}b(n)\, n^{-s}
=
(1+q^{-s})\, \zeta(2s)\, \sum_{n\ge 1}\, |a_f(n)|^2 n^{-s}
\leqno(15)
$$
is entire except for a simple pole at $s=1$ with 
$$
\mathop{{\rm Res}}_{s=1}\, \varphi(s) = \frac{2\pi^2}{q}\, \langle f, f\rangle
$$
and satisfies the functional equation
$$
\Phi(s) = 
\left(\frac{q}{4\pi^2}\right)^s\, \Gamma(s)^2 \varphi(s) = \Phi(1-s)
\leqno(16)
$$
by \cite{L}. Suppose that $F \in C^{\infty}_c(0, \infty)$ has
$\int_0^\infty F(x) \, dx= 1$. Then the Mellin transform
$$
\hat F(s) = \int_0^\infty F(x) \, x^s \, \frac{dx}{x}
$$
is entire, of rapid decay in vertical strips and  $\hat F(1)= 1$.
By Mellin inversion,
$$
F(x)= \frac{1}{2\pi i} \int\limits_{{\rm Re}\, s=2} \hat F(s) x^{-s}ds.
$$
for $x>0$. Thus 
$$
\eqalign{
\sum_{n\ge 1} b(n) F(n/q^2) 
&=
\frac{1}{2\pi i } \int\limits_{{\rm Re}\, s=2} \hat F(s)\varphi(s)
q^{2s} ds
\cr\noalign{\vskip 5pt}
&=
2\pi^2 q\, \langle f, f\rangle + \frac{1}{2\pi i}
\int\limits_{{\rm Re}\, s = -1} \hat F(s) \varphi(s) q^{2s} ds.
}
\leqno(17)
$$
It follows from (16), (15) and (12) that 
$$
|\varphi(-1+it)| \leq (4\pi^2)^{-3} q^3
\frac{|\Gamma(2-it)|^2}{|\Gamma(-1+it)|^2} (1+q^{-2})|\zeta(4-2it)|
\sum_{n \geq 1} d(n)^2 n^{-2}
$$
and by Stirling's formula this is 
$ \ll q^3 (|t|+ 1)^6$
with an absolute constant.  Hence
$$
\frac{1}{2\pi i}\int_{{\rm Re}\,(s)= -1}\hat{F}(s)\varphi(s)q^{2s}ds
\ll q \int_{-\infty}^{\infty}|\hat{F}(-1+it)|(|t|+1)^6dt \ll q
\leqno(18)
$$
since $\hat{F}(-1+it)\ll (1+|t|)^{-8}$.
By (17) and (18) we get
$$
\langle f, f\rangle = \frac{1}{2\pi^2 q}\, \sum_{n\ge 1} b(n)
\, F(n/q^2) + O(1)\,
$$
where the implied constant depends only on $F$. Now, by (15) and (12)
$$
0\le b(n) \le 2\, \sum_{n=m\ell^2} d^2(m) = 2d_4(n)
$$
where $d_4(n)$ is the number of factorizations of $n$ into $4$ factors. Thus
$$
\langle f,f\rangle 
\ll 
\frac{1}{q}\, \sum_{n\ge 1} d_4(n)\, F(n/q^2)\ll q\, \log^3 q.
$$
Fixing $F$ gives the result with an absolute implied constant. $\Box$

\bigskip\noindent
(b) The fact that $f$ cannot be tetrahedral is shown in 
\cite[Theorem 7(c)]{S}. It is also observed in \cite[p.~362]{S}
that if $p\ne q$ and $f\in N_{{\rm oct}}$ then
$$
\left(\frac{p}{q}\right)\, a(p^2)\in \{-1, 0, 1, 3\}\,,
\leqno(19)
$$
while if $f\in N _{{\rm ico}}$ then 
$$
\left(\frac{p}{q}\right)\, a(p^2)
\in 
\left\{-1, 0, 3, \frac{1+\sqrt{5}}{2},\frac{1-\sqrt{5}}{2}\right\}\,.
\leqno(20)
$$
For general $f\in N_1(q, \e)$ it follows from (11) that if we set
$x=\bar \e(p) \, a_f(p^2)$ then for $p\ne q$ and $n\ge 0$
$$
\bar \e ^n (p)\, a_f(p^{2n})=P_n(x)\,,
\leqno(21)
$$
where $P_0(x)=1$, $P_1(x)=x$, and 
$P_{n+1}(x) = (x-1)\, P_n(x) - P_{n-1}(x)$ for $n\ge 1$. 
Thus $P_2(x) = x^2 - x - 1$, $P_3(x) = x^3 - 2 x^2 - x + 1$, etc.
One may check that
$$
x(x+1)(x-1) (x-3) + 1
=
P_4(x) - P_2(x) - P_1(x)\,,
$$
while
$$
x(x+1)(x-2) (x-3) (x^2-x-1) + 1
=
P_6(x) - P_4(x) -P_1(x)\,.
$$
Thus by (19)--(21) we finish the proof of (b).~$\Box$

\medskip\noindent
{\bf Remark:} It has been checked that no linear form $\sum_{\ell=1}^7
c_\ell\, a_f(p^{\ell})$ takes positive values for all possible values
of $a_f(p)$, $p\ne q$, when $f\in N_{{\rm oct}}$.  A similar remark
applies to $f\in N_{{\rm ico}}$.

\subsection*{4. Counting newforms and quartic fields}

By combining Propositions 1 and 2 we may now estimate the number of
newforms of octahedral or icosahedral type.

\medskip\noindent
{\bf Proposition 3.} {\sl For $q$ prime 
$$
\leqalignno{
&\# N_{{\rm oct}} \ll q^{7/8} \log^4 q
&{\rm (a)}
\cr\noalign{\vskip 5pt}
&\# N_{{\rm ico}} \ll q^{11/12} \log^4 q
&{\rm (b)}
\cr
}
$$
with absolute implied constants.
}

\medskip\noindent
{\bf Proof:} By Propositions 1 and 2(a) we have for any $c_n \in {\bf C}$
$$
\sum_{f\in N_1(q)}\, |\sum_{n\le N}c_n\, a_f(n)|^2
\ll
(q+N)\, \log^3 q \sum_{n\le N}|c_n|^2\,.
$$
Thus we deduce by positivity the inequality 
$$
\sum_{f\in N_{{\rm oct}}}\, |\sum_{n\le N}c_n\, a_f(n)|^2
\ll
(q+N)\, \log^3 q \sum_{n\le N}|c_n|^2.
\leqno(22)
$$
Choose $N=q$ and 
$$
c_n =\cases{ 1, & $n=p^8\le q$,\cr
             -1, & $n=p^4\le q^{1/2}$,\cr
        -\left(\frac{p}{q}\right), & $n=p^2 \le q^{1/4}$, \cr
          0, & {\rm otherwise}, \cr
}
$$
for $p$ prime. By the prime number theorem and (13)
$$
\sum_{n\le N} c_n a_f(n)
\sim
\frac{8q^{1/8}}{\log q},\qquad {\rm for }\; f\in N_{{\rm oct}}\,,
$$
while 
$$
\sum_{n\le N}|c_n|^2 \sim \frac{24 q^{1/8}}{\log q}\,.
$$
Hence, by (22)
$$
\# N_{{\rm oct}}\ll q^{7/8}\, \log^4 q\,.
$$
A similar argument works for $N_{{\rm ico}}$.~$\Box$

\medskip
Theorem 1 now follows from Propositions 2(b) and 3 since
$$
\# N_{{\rm dih}} = \frac{h-1}{2}\ll q^{1/2}\log  q\,.
$$

\medskip
For its Corollary we use the fact, shown in \cite{S}, that the 
required space of differentials has dimension
$$
\frac{1}{2}\, {\rm dim}\, S_1(q) - \frac{1}{4}\, (h-1)\,.
$$

\medskip
Theorem 2 follows from Proposition 3(a), Theorem 8 of \cite{S} and
Tunnell's proof \cite{T} of the Artin conjecture for octahedral Galois
representations. Here we also use the fact that the Galois closure of
any quartic extension with discriminant $-q$ has Galois group $S_4$,
since otherwise the discriminant would not be square free (e.g., see
\cite{B}).

\medskip\noindent
{\it Acknowledgement.} I would like to thank J-P. Serre for his
many helpful comments.  This paper was written while I was participating in
the 1994-95 program on Automorphic Forms at MSRI.  I thank the
organizers and the Institute for providing financial support and a
congenial atmosphere in which to work.


\bigskip\bigskip
\renewcommand{\baselinestretch}{0.5}
\noindent
{\sc William Duke}

Department of Mathematics, Rutgers University 

New Brunswick, NJ 08903 

{\it email\/}: duke\verb+@+math.rutgers.edu
\bigskip

\noindent{\it Current Address\/}:

Mathematical Sciences Research Institute

1000 Centennial Drive 

Berkeley, CA 94720--5070

\end{document}